\documentclass[12pt,a4paper]{amsart}
\usepackage{graphicx} 
\usepackage{caption}

\usepackage{algorithm}
\usepackage{algpseudocode}\usepackage{xcolor}
\usepackage{subcaption}
\usepackage{newfloat}
\usepackage{listings}
\usepackage{amsmath}
\usepackage{amsfonts}
\usepackage{amsthm}

\usepackage{multirow}
\usepackage[colorlinks,pdfdisplaydoctitle,linkcolor=blue,citecolor=blue]{hyperref}
  
\usepackage[margin=1.2in]{geometry}
\usepackage{comment}

\newcommand{\bx}{\boldsymbol{x}}
\newcommand{\by}{\boldsymbol{y}}
\newcommand{\bu}{\boldsymbol{u}}
\newcommand{\ba}{\boldsymbol{a}}
\newcommand{\bI}{\boldsymbol{\mathbf{I}}}
\newcommand{\bZ}{\boldsymbol{0}}
\newcommand{\beps}{\boldsymbol{\epsilon}}
\newcommand{\bmu}{\boldsymbol{\mu}}

\newcommand{\bs}{\boldsymbol{s}}

\newcommand{\bG}{\boldsymbol{\mathcal{G}}}

\newcommand{\mL}{\mathcal{L}}
\newcommand{\mN}{\mathcal{N}}

\newcommand{\mI}{\mathcal{I}}
\newcommand{\mS}{\mathcal{S}}

\newcommand{\RN}[1]{\textup{\uppercase\expandafter{\romannumeral#1}}
}

\theoremstyle{definition}

\begin{document}

\title[Generative Downsampling with PGDM]{Generative Downscaling of PDE Solvers With  Physics-Guided Diffusion Models}

\author{Yulong Lu}

\address{(YL) School of Mathematics, University of Minnesota Twin Cities, Minneapolis, MN 55414.}

\email{yulonglu@umn.edu}

\author{Wuzhe Xu} 

\address{(WX) Department of Mathematics and Statistics, University of Massachusetts Amherst, Amherst, MA 01003.}

\email{wuzhexu@umass.edu}



\begin{abstract}
Solving partial differential equations (PDEs) on fine spatio-temporal scales for high-fidelity solutions is critical for numerous scientific breakthroughs. Yet, this process can be prohibitively expensive, owing to the inherent complexities of the problems, including nonlinearity and multiscale phenomena. To speed up large-scale computations, a process known as downscaling is employed, which generates high-fidelity approximate solutions from their low-fidelity counterparts. In this paper, we propose  a novel  Physics-Guided Diffusion Model (PGDM) for downscaling. Our model, initially trained on a dataset comprising low-and-high-fidelity paired solutions across coarse and fine scales, generates new high-fidelity approximations from any new low-fidelity inputs. These outputs are subsequently refined through fine-tuning, aimed at minimizing the physical discrepancies as defined by the discretized PDEs at the finer scale. We evaluate and benchmark our model's performance against other downscaling baselines in three categories of nonlinear PDEs. Our numerical experiments demonstrate that our model not only outperforms the baselines but also achieves a computational acceleration exceeding tenfold, while maintaining the same level of accuracy as the conventional fine-scale solvers.  
\end{abstract}

\maketitle

\section{Introduction}
Numerical simulation of PDEs play an indispensable role in science and engineering. Traditional numerical methods, such as finite difference method and finite element method, often become computationally intensive with an increase in mesh grids. This increase is typically necessary to accurately resolve PDEs, given their complexities, such as nonlinearity, scale-separation stiffness, and high dimensionality. In recent times, the adoption of deep learning techniques to develop more efficient numerical methods has seen a significant rise in popularity. Numerous studies have explored the direct approximation of solutions using neural networks. The work \cite{raissi2019physics} proposed the physics-informed neural networks (PINNs) that minimizes the $L^2$-loss associated with the governing physics, and it has proven to be highly efficient in addressing various complex PDE problems, such as the fluid dynamic \cite{cai2021physics}, inverse problem \cite{yu2022gradient} and multiscale problem \cite{lu2022solving, jin2023asymptotic}.  Furthermore, it is worth mentioning several variations and enhancements of PINNs, such as \cite{kharazmi2019variational, kharazmi2021hp, wang2023expert}. 
Many other approaches have also been developed, such as Deep Ritz Method \cite{weinan2018deep,lu2021priori,lu2022priori}, based on the variational (or Ritz) formulation of PDEs, the deep BSDE method \cite{han2017deep,han2018solving} for certain class of parabolic PDEs based on their probabilistic and control formulation, and the weak adversarial networks \cite{zang2020weak} based on the Galerkin (or weak) formulation. Additionally, deep learning has been leveraged to expedite classical iterative solvers \cite{chen2022meta,um2020solver,hsieh2018learning,arisaka2023principled,azulay2022multigrid,nikolopoulos2024ai}, showcasing its versatility and potential in enhancing computational efficiency across various PDE-solving methodologies.

In many scientific domains, there's a notable interest in discerning mappings or operators between infinite-dimensional function spaces. Recent advancements have seen neural networks being harnessed to approximate these solution operators. Among the leading neural operator models are the Deep operator networks (DeepONet) \cite{lu2019deeponet, wang2023long} and the Fourier neural operator (FNO) \cite{li2020fourier, goswami2022physics, li2021physics}. Nevertheless, these standard neural operators often necessitate a substantial training dataset composed of numerous parameter-solution pairs, posing challenges in scenarios where solution labels are costly to obtain. To circumvent this issue, physics-informed neural operators \cite{wang2021learning,goswami2023physics}  have been introduced, merging PDE constraints with operator learning by embedding known differential equations directly into the training loss function.
Notably, neural operators have found applications in downscaling climate data \cite{jiang2023efficient,yang2023fourier}  and enhancing super-resolution in imaging \cite{wei2023super}, as well as general PDE problems \cite{kovachki2021neural}, showcasing their versatility.
  The method we propose in this paper will also be compared with the neural operator baselines.

In this paper, we focus on accelerating the computation of PDEs from a downscaling viewpoint. In climate modeling and simulation, downscaling \cite{wilby1998statistical} refers to a class of methods that generate high-fidelity climate data out of their low-fidelity counterpart. Similar processes may carry with different names. For instance in the community of imaging and computer vision, this process is named super-resolution. Such a downscaling/super-resolution is appealing because low-fidelity solutions can be generated via solving PDEs on coarse-grids which is computationally much cheaper compared to the high-fidelity solutions.
Classical downscaling techniques in climate science and meteorology have ranged from pointwise regression \cite{sachindra2018statistical,vandal2019intercomparison} to super-resolution \cite{vandal2017deepsd} and maximum likelihood estimation \cite{bano2020configuration}. Recent initiatives  have seen Fourier neural operators being applied for downscaling \cite{yang2023fourier}, effectively bridging fast low-resolution simulations to high-resolution climate outputs. Moreover, the use of deep generative models for climate data downscaling, inspired by their success in computer vision for super-resolution, has gained traction.
In   \cite{leinonen2020stochastic,price2022increasing},  Generative Adversarial Networks (GANs) were adopted for downscale precipitation forecasting. In \cite{groenke2020climalign}, the authors proposed a ClimAlign approach to downscaling with normalizing flows. Recently, 
diffusion models have demonstrated their ability to produce high quality samples,    beating many competing generative models such as GANs in numerous machine learning problems \cite{dhariwal2021diffusion}. Specifically,  diffusion-models \cite{ho2020denoising,song2020score,song2020denoising} are capable of generating high-fidelity (super-resolution) images \cite{ho2022cascaded}. Inspired by their great success in machine learning tasks, we propose to deploy diffusion models for downscaling PDE solvers. Unlike the purely data-driven nature of downscaling in the aforementioned applications, our PDE-focused generative downscaling necessitates adherence to the physical laws governing the PDE model, integrating a unique challenge to this innovative approach.

\subsection{Related work}

The utilization of deep generative models for solving PDEs is not new. Initial studies, such as those by \cite{farimani2017deep,joshi2019generative}, harnessed adversarial generative models to tackle PDEs by integrating a physics-informed loss with a uniquely tailored adversarial loss. More recent efforts, like the study by \cite{yang2023denoising}, showcased the diffusion model's efficacy as an alternative solution operator, mapping  initial conditions to solutions at subsequent times, and displaying competitive prowess alongside other neural operator approaches. To overcome the hurdle of scarce data in operator learning, \cite{apte2023diffusion}  utilized the diffusion model for the creation of synthetic data samples, thereby enriching the training dataset.
 The research work \cite{shu2023physics} is most relevant to ours, where the authors developed a physics-informed diffusion model designed to accurately reconstruct high-fidelity samples from low-fidelity ones. While both our approach and the method proposed in \cite{shu2023physics} intend to reconstruct high-fidelity data from low-fidelity sources, our approach differs from \cite{shu2023physics} in several aspects. First, they added physics-informed loss in the denosing score-matching loss in the training of their diffusion models, which can be very expensive due to the enforcement of PDE information in each gradient step during the training. In contrast, our physics-guided diffusion model decouples the step of purely data-driven conditional diffusion model from the physics-enhancement step. The pre-trained model when combined with low-fidelity input produces a high-fidelity output that can be used as a good warm-start for minimizing the physics-informed loss in the second step.  This two-step procedure improves substantially the efficiency of  training and accuracy of generated solutions. Source code is available at \href{https://github.com/woodssss/Generative-downsscaling-PDE-solvers}{https://github.com/woodssss/Generative-downsscaling-PDE-solvers}

\subsection{Our contributions.} We introduce the physics-guided diffusion model as a universal framework for downscaling PDE solutions from low-resolution to high-resolution. 
\begin{itemize}

    \item We first reformulate the downscaling problem as a conditional sampling task, where the objective is to sample from the posterior distribution of unknown high-fidelity solutions, given any arbitrary low-fidelity input. This reformulation allows for a more targeted and accurate generation of high-resolution outputs from their lower-resolution counterparts.

    \item  The first step involves conditional sampling via a diffusion model to produce preliminary high-fidelity samples. Subsequently, these samples are refined through a physics-informed loss minimization step, ensuring they adhere to the physical laws governing the PDEs. This dual-step approach effectively merges data-driven sample generation with physics-based accuracy enhancement.

    \item The proposed method consistently outperforms several existing downscaling baselines in a range of nonlinear static and time-dependent PDEs. Remarkably, it not only matches the accuracy of traditional high-fidelity solvers at the fine scale but also achieves this with a significant reduction in computational expenditure, cutting costs by more than tenfold.

    
    
\end{itemize}

\section{Problem Set-Up} 
\subsection{Problem description}
 Our primary focus is on developing efficient approximations for the solution 
$u$ to a generic PDE, subject to appropriate boundary conditions, as outlined below:
\begin{equation}\label{eqn:pde}
    \mL u = a,
\end{equation}
where $\mL$ is the (possibly nonlinear) differential operator  and $a$ is the source term. Traditional grid-based PDE solvers typically approach this by solving a discretized version of the problem:
$$
\mathcal{L}_h \bu_h = \ba_h + \boldsymbol\epsilon_h.
$$
Here $h$ represents the spatial-temporal grid-size, $\boldsymbol\epsilon_h$ denotes certain (unknown) noise that potentially encapsulating errors in the pointwise evaluations of the functions, and $\mL_h$,  $\bu_h$ and $\ba_h$ represent the discrete  approximations to the operator $\mL$, the solution $u$ and the source $a$ respectively. Discretization on fine grids (characterized by a small 
$h$) usually results in high-fidelity (high-resolution) solutions but at the expense of significantly increased computational costs. Therefore, finding an optimal cost-accuracy balance is crucial. Downscaling, in this context, refers to a series of techniques that first solve PDEs on coarse grids and subsequently convert the low-fidelity solutions obtained on these coarse grids to their high-fidelity equivalents on fine grids, offering a strategic approach to manage the trade-offs between computational expense and solution accuracy.

\subsection{Downscaling as conditional sampling}To describe our diffusion-based downscaling approach, we would like to first present a conditional-sampling formulation to the downscaling problem. To ease the notation, we will suppress the dependence of quantities on the grid size $h$ and  denote by $\bu^c$ and $\bu^f$  the low-fidelity solution and high-fidelity solution respectively. Similarly one can define for $g = c, f$ the operators $\mL^g$, the source terms $\ba^g$ and the noise $\boldsymbol\epsilon^g$. Moreover,  we have 
$$
\mathcal{L}^g \bu^g = \ba^g + \boldsymbol\epsilon^g.
$$
Assume that $a^c = \mathcal{R} a^f$ with some fine-to-coarse restriction operator $\mathcal{R}$, one has 
$$
\mathcal{R} \mL^f \bu^f = \mL^c \bu^c + \mathcal{R} \boldsymbol\epsilon^f - \boldsymbol\epsilon^c. 
$$
Assuming the invertibility of $\mL^c$, we can rewrite above as  
$$
\bu^c =  (\mL^c)^{-1}\mathcal{R} \mL^f \bu^f +  \boldsymbol\epsilon, \text{ where } \boldsymbol\epsilon := -(\mL^c)^{-1} ( \mathcal{R} \boldsymbol\epsilon^f - \boldsymbol\epsilon^c).
$$
In another words, the downscaling problem is an  inverse problem of recovering $\bu^f$ from the noisy downscaling observation $\bu^c$ via 
$$
\bu^c = \bG \bu^f +\boldsymbol\epsilon, \text{ where } \bG:=(\mL^c)^{-1}\mathcal{R} \mL^f . 
$$
We adopt the Bayesian approach for solving the inverse problem.   Given a prior $p(\bu^f)$ on the set of fine solutions, one can define by the Bayes' rule the posterior distribution 
\begin{equation}\label{eq:post}
 p(\bu^f | \bu^c) \propto p ( \bu^c | \bu^f)\times p(\bu^f),
\end{equation}
where $p( \bu^c | \bu^f)$ is the likelihood function. With above, we have recast the downscaling as the problem of conditional sampling from the Bayesian posterior $p(\bu^f | \bu^c)$ given an arbitrary low-fidelity input $\bu^c$. 

\subsection{Challenges and our approach}\label{sec:challenge} Despite the appealing conditional sampling framework offered by downscaling, direct sampling from the Bayesian posterior \eqref{eq:post} presents infeasibility and numerous challenges.  First, the prior $p(\bu^f)$ is unknown and needs to be learned from the data.
Secondly, evaluating the likelihood function poses significant difficulties due to two primary reasons: (1) the forward map $\bG$ is either inaccessible, owing to an unknown fine-to-coarse restriction or, even if known, the computation may be prohibitively expensive, and (2) the noise distribution is typically unknown, leading to an intractable likelihood. Last but not the least, the conditional samples $\bu^f$ given $\bu^c$ even can be generated may not fulfil the discrete PDE problem on the fine grid, especially given a limited amount of data.  To address these issues, we introduce a physics-guided diffusion model designed to learn and draw physics-conformal high-fidelity samples from any low-fidelity inputs.
Our approach is methodically divided into two pivotal steps:

\begin{itemize}
    \item[(1)] {\bf Pre-training step}: We pre-train a conditional diffusion model using a finite collection of low-and-high fidelity solution pairs $\{(\bu^c_i, \bu^f_i)\}_{i=1}^n$  laying the groundwork for subsequent refinements. 

    \item[(2)] {\bf Refining step}: Upon receiving any new low-fidelity input, we refine the output via the pre-trained model to ensure an enhanced fit with the fine-grid PDE, thereby further improving the solution's fidelity.
\end{itemize}

The pre-training phase of our approach is primarily data-driven and accounts for the majority of computational expenditure.  In contrast, the refining step is more computationally economical and aims to enhance the high-fidelity output by minimizing the physics misfit loss. This enhancement could be achieved, for instance, by executing few, such as two, Gauss-Newton iterations, starting with the initial output from the pre-trained model, thereby streamlining the process towards achieving superior solution accuracy with reduced computational demand.

\section{Methodology}
\subsection{Unconditioned Diffusion model}
To introduce our conditional diffusion models for downscaling, we start with a general overview of unconditioned diffusion models, with our focus on the Denoising Diffusion Probabilistic Models and one of its accelerated version called Denoising Diffusion Implicit Models.

\subsubsection{Denoising Diffusion Probabilistic Models(DDPM)} Let $q(\bx)$  represent the target data distribution.  DDPM constructs a Markovian noising process that incrementally contaminates the data $\bx_0$ with Gaussian noise over 
$T$ steps, ultimately transforming it into pure Gaussian noise. This noising process is denoted  by $q(\bx_{0:T})$, where the $\bx_1 , \cdots, \bx_T$ are progressively noised versions of the data, all maintaining the same dimensionality as $\bx_0$, and $q(\bx_T)$ is approximately an isotropic Gaussian distribution. This forward process of $\bx_{0:T}$ can be described by the Markov process with the transition kernel defined by 
\begin{equation}
    q(\bx_t | \bx_{t-1}) := \mN(\bx_t; \sqrt{\alpha_t} \bx_{t-1}, (1-\alpha_t) \bI)
\end{equation}
where $\{ \alpha_t \}_{t=0}^T \subset (0,1)$ is a sequence of user designed parameters. Since the noises we add in each step are Gaussian, we have
\begin{equation}\label{eqn:forw_t_0}
\bx_t | \bx_{0} \overset{d}{=} \sqrt{\bar{\alpha}_t} \bx_0 + \beps \sqrt{1-\bar{\alpha}_t}, \quad \beps \sim \mN(0, \bI),
\end{equation}
where $\bar{\alpha}_t = \Pi_{t=1}\alpha_t$. Generation of new data samples can be done via the backward
process in DDPM. More precisely, with the assumption that the reverse process $q(\bx_{t-1}|\bx_t)$ can be modeled as Gaussians with trainable mean and fixed variance, a reversed Markov process is parameterized in
the form of
\begin{equation}\label{eqn:trans_kernal}
    p_{\theta}(\bx_{t-1}|\bx_t) := \mN \Big(\bx_{t-1}; \frac{1}{\sqrt{\alpha_t}}\Big(\bx_t + \frac{\sqrt{(1-\alpha_t)}}{\sqrt{(1-\bar{\alpha}_t)}}\bs_\theta (\bx_t, t)\Big), \sigma_t^2 \bI\Big),
\end{equation}
where $\sigma_t^2 =\frac{(1-\alpha_t)(1-\bar{\alpha}_{t-1})}{1-\bar{\alpha}_{t}}$
and is trained with the weighted  evidence lower bound (ELBO) 
$$
\hat{\theta} = \arg\min_{\theta} \sum_{t=1}^T (1-\bar{\alpha}_t)\mathbf{E}_{q(\bx_0)} \mathbf{E}_{q(\bx_t | \bx_0)} \Big\| \frac{1}{\sqrt{1-\bar{\alpha}_t}}\bs_\theta (\bx_t, t) - \nabla_x \log(p(\bx_t|\bx_0))\Big\|^2. 
$$
It can be shown further by integration by parts that minimizing the ELBO is
equivalent to the denoising problem
$$
\hat{\theta} = \arg\min_{\theta} \sum_{t=1}^T \mathbf{E}_{\bx_0\sim q(\bx_0)} \mathbf{E}_{\beps_t \sim \mathcal{N}(0,\mathbf{I}) } \Big\|  \frac{(1-\alpha_t)^2}{2\sigma_t^2 \alpha_t (1- \bar{\alpha}_t)}
 \bs_\theta (\sqrt{\bar{\alpha}_t} \bx_0 + \beps_t \sqrt{1-\bar{\alpha}_t}, t) + \beps_t\Big\|^2. 
$$
The optimized neural network $\bs_{\hat{\theta}}$ enables us to generate new samples $\bx_0$ through  the backward process iterates: starting with $\bx_T\sim \mathcal{N}(0, \mathbf{I})$, 
$$
\bx_{t-1} = \frac{1}{\sqrt{\alpha_t}}\Big(\bx_t + \frac{\sqrt{(1-\alpha_t)}}{\sqrt{(1-\bar{\alpha}_t)}}\bs_\theta (\bx_t, t) \Big) +  \sigma_t \xi_t, \quad t=T, T-1, \cdots, 1,
$$
where $\{\xi_t\}_{t=1}^T \overset{i.i.d.}{\sim} \mathcal{N}(0, \mathbf{I})$. 



\subsubsection{Accelerating sampling with Denoising Diffusion Implicit Models (DDIM)} One major drawback of  DDPM is that generating a new sample from the data distribution requires  simulating the whole Markov backward process  for many (typically hundred or thousand) steps (or equivalently the forward network passes), which can be computationally intensive and time-consuming. Recently, Song et. al. \cite{song2020denoising} proposed the denoising diffusion implicit models (DDIM), accelerating the generative process by using a non-Markovian deterministic diffusion pathway, culminating in implicit models capable of producing samples at an faster pace without compromising on quality. More specifically, given a selective increasing sequence of length $L$, denoted by  $\{\tau_i\}_{i=1}^L \subset [1, 2, \cdots, T]$, DDIM generates a sample $\bx_{\tau_{i-1}}$ from $\bx_{\tau_i}$ by making the following update:



\[
\displaystyle \bx_{\tau_{i-1}} = \frac{\sqrt{\bar{\alpha}_{\tau_{i-1}}}}{\sqrt{\bar{\alpha}_{\tau_{i}}}} \bx_{\tau_{i}} + (\frac{\sqrt{1-\bar{\alpha}_{\tau_{i}}} \sqrt{\bar{\alpha}_{\tau_{i-1}}}}{\sqrt{\bar{\alpha}_{\tau_{i}}}} 
 - \sqrt{1 - \bar{\alpha}_{\tau_{i-1}}})\bs_{\hat{\theta}}\left(\mathbf{x}_{\tau_{i}}, {\tau_{i}} \right), \quad i = 1, 2, \cdots, L, 
\]
where $\bs_{\hat{\theta}}$ is the optimal score network trained in the same manner as in DDPM.

\subsection{Conditioned diffusion model}
Now let us move on to the problem of conditional sampling using conditional diffusion model. Recall that our goal is to sample from the posterior distribution 
$$ 
p(\bu^f | \bu^c) \propto p ( \bu^c | \bu^f)\times p(\bu^f)
$$
for any given low-fidelity solution $\bu^c$. One straightforward idea for doing so would be to train a conditional score network $\bs_{\theta}(\bu^f, \bu^c, t)$ that minimizes the ELBO in the conditional setting: 
$$
\hat{\theta} = \arg\min_{\theta} \sum_{t=1}^T (1-\bar{\alpha}_t)\mathbf{E}_{\bu_0^f\sim p(\bu^f)} \mathbf{E}_{p(\bu_t^f | \bu_0^f)} \Big\|\frac{1}{\sqrt{1-\bar{\alpha}_t}}\bs_\theta (\bu_t^f, \bu^c, t) - \nabla_{\bu_t^f} \log(p(\bu_t^f|\bu_0^f, \bu^c))\Big\|^2. 
$$
By the Bayes' formula, the true conditional score function $\nabla_{\bu^f} \log(p(\bu^f|\bu_0^f, \bu^c))$ can be written as 
$$
\nabla_{\bu^f} \log(p(\bu^f|\bu_0^f, \bu^c)) = \nabla_{\bu^f} \log p(\bu^f) + \nabla_{\bu^f} \log(p( \bu^c|\bu_0^f,\bu^f)),
$$
where the first term represents the score function corresponding to the prior $p(\bu^f)$ and the second term encodes the conditional likelihood. While the prior $p(\bu^f)$ can be learned from high-fidelity training samples, the conditional likelihood is often computationally intractable and existing conditional diffusion models resort to various approximations to the conditional likelihood, such as the pseudo-inverse in the inverse problem setting or the posterior mean in the general nonlinear inverse problem setting. Unfortunately, it is  impossible to construct those approximations in our setting 
due to the lack of the complete knowledge on the forward operator $\bG$ as we illustrated in Section \ref{sec:challenge}. To bypass these issues, we seek a purely data-driven approach to learn the conditional score without incorporating the forward model. To be concrete, given a training set of low-and-high fidelity solution pairs $\{(\bu^c_k, \bu^f_k)\}_{k=1}^N$, we seek to optimize  the score network $\bs_\theta (\bu^f, \bu^c, t)$ with respect to the parameter $\theta$ such that
\begin{equation}\label{eq:condsgm}
    \hat{\theta} = \arg\min_{\theta} \sum_{t=1}^T  \frac{1}{N}\sum_{k=1}^N  \Big[ \frac{(1-\alpha_t)^2}{2\sigma_t^2 \alpha_t (1- \bar{\alpha}_t)} \Big\|\bs_\theta (\sqrt{\bar{\alpha}_t} \bu_k^f + \beps_{t,k} \sqrt{1-\bar{\alpha}_t},\bu_k^c, t) + \beps_{t,k} \Big\|^2 \Big],
\end{equation}
where  $\{\beps_{t,k}\} \overset{i.i.d.}{\sim} \mathcal{N}(0, \mathbf{I}), t=1,\cdots, T; k=1, \cdots, N$. As discussed in \cite{ho2020denoising}, it is  beneficial to sample quality and simpler to implement to omit the time dependent coefficient $\frac{(1-\alpha_t)^2}{2\sigma_t^2 \alpha_t (1- \bar{\alpha}_t)}$, and the training process for the the conditional diffusion model is summarized in Algorithm~\ref{alg:train} below.

Similar to the unconditioned setting, with the optimal score network $\bs_\theta (\bu_t^f, \bu^c, t)$, we can generate a new high-fidelity sample $\bu^f= \bu^f_0$ conditioned on a new low-fidelity solution $\bu^c$ by evolving the backward process with a terminal sample $\bu^f_T \sim \mathcal{N}(0, \mathbf{I})$. In the framework of DDPM, such a backward process is given by  
$$
\bu^f_{t-1} = \frac{1}{\sqrt{\alpha_t}}\Big(\bu^f_t + \frac{\sqrt{(1-\alpha_t)}}{\sqrt{(1-\bar{\alpha}_t)}} \bs_{\hat{\theta}} (\bu^f_t,\bu^c, t)\Big) +  \sigma_t \xi_t, \quad t=T, T-1, \cdots, 1.
$$
where $\{\xi_t\}_{t=1}^T \overset{i.i.d.}{\sim} \mathcal{N}(0, \mathbf{I})$.

In the case of DDIM, the backward process updates according to 
$$
\displaystyle \bu^f_{\tau_{i-1}} = \frac{\sqrt{\bar{\alpha}_{\tau_{i-1}}}}{\sqrt{\bar{\alpha}_{\tau_{i}}}} \bu^f_{\tau_{i}} + (\frac{\sqrt{1-\bar{\alpha}_{\tau_{i}}} \sqrt{\bar{\alpha}_{\tau_{i-1}}}}{\sqrt{\bar{\alpha}_{\tau_{i}}}} 
 - \sqrt{1 - \bar{\alpha}_{\tau_{i-1}}})\bs_{\hat{\theta}}\left(\bu^f_{\tau_{i}}, \bu^c,  {\tau_{i}} \right), \quad i = 1, 2, \cdots, L.
$$
In practice, we observe that incorporating the information of the source term $a$ in the training of conditional score network can improve the sample quality in the sense of better fitting the PDE on the fine scale. Therefore throughout the paper we look for a score network depend on $\ba^f$ that solves \eqref{eq:condsgm} with the score network $\bs_\theta (\bu^f, \bu^c, t)$ replaced by $\bs_\theta (\bu^f, \bu^c, \ba, t)$.

\begin{algorithm}
	\caption{Training of conditional diffusion models} 
    \label{alg:train}
	\begin{algorithmic}[1]
 \Require Training dataset $\mathcal{S}:=\{ 
\bu^c_k, \bu^f_k, \ba_k \}_{k=1}^{N}$, hyperparameter $\{\alpha_t\}_{t=0}^T \subset (0,1)$, batch size $B$.
        \Repeat
        \State Sample $\{\bu_j^c, \bu_j^f, \ba_j \}_{j=1}^B \sim S$, let $\bu^f_{0, j} = \bu_j^f, ~ j = 1, \cdots, B$
        \State $t \sim $ Uniform$(\{ 1, \cdots, T \})$
        \State $\beps_{t,j} \sim \mN(\bZ, \bI), ~ j = 1, \cdots B$
        \State Compute $\bu^f_{t, j} = \sqrt{\bar{\alpha}_t} \bu^f_{0, j} + \sqrt{1-\bar{\alpha}_t} \beps_j$
        \State Take gradient descent step on \\
        \quad \quad \quad \quad \quad $\nabla_{\theta} \Big[ \sum_{t=1}^T  \frac{1}{B}\sum_{j=1}^B \|\bs_\theta (\bu^f_{t, j},\bu_{i}^c, t) + \beps_{t,j} \|^2 \Big]$
    \Until{converged.}
	\end{algorithmic} 
\end{algorithm}

\subsubsection{Physics-guided diffusion model}
The generated high-fidelity sampled solution through the conditional diffusion model contains rich information from the training high-fidelity training data while informed by the low-fidelity input. Yet, the generated sample may not fulfil the PDE on the fine scale and hence need to be further enhanced to better conform with the physics. To improve the solution quality, we refine the solution by solving the (nonlinear) least square problem
\begin{equation}\label{eqn:min}
    \bu \in \arg \min_{\bu} \| \mL^f \bu - \ba \|_2,
\end{equation}
where the boundary term in the loss for simplicity and could  be included in practice. In our experiments, we generate a refined solution by solve problem \eqref{eqn:min} with few Gauss-Newton iterations and initial guess chosen as the generated output from the conditioned diffusion model. 

By combining the conditional sample generation step with condition diffusion models with the refining step with Gauss-Newton, we present below the overall physics-guided diffusion model for downscaling. Assume that we have access to a pre-trained conditional diffusion model or specifically the conditional score function $\bs_\theta$ (see Algorithm 1). We present the algorithm~\ref{alg:fast_sample} for the refined sample generation process in the framework of DDIM. Through our experiments, we have determined that a refining step of $t_s = 1$ sufficiently enhances the solution quality to be on par with that of the fine solver. For clarity, the $\bu_{t,j}^f$, in Algorithm~\ref{alg:train}, denotes the $j$th noised data at diffusion time step $t$, while the $\bu^f_{\tau_{i-1}}$ in Algorithm~\ref{alg:fast_sample} denotes the reconstruction at diffusion time step $\tau_{i-1}$. Furthermore, we employ the Gaussian-Newton algorithm to refine solutions produced by the diffusion model. Detailed implementation can be found in Algorithm~\ref{alg:gn} in the Appendix.
\begin{algorithm}[H]
	\caption{Physics-guided diffusion model (PGDM) for downscaling} 
    \label{alg:fast_sample}
	\begin{algorithmic}[1]
    \Require A given low-fidelity $\bu^c$ and the source $\ba^f$ evaluated on the fine scale, hyperparameters $\{\alpha_t\}_{t=0}^T \subset (0,1)$ and a set of indices $\{\tau_i\}_{i=1}^L \subset [1, 2, \cdots, T]$ with length $L$.
    \State $\bu^f_T \sim \mN(\bZ, \bI)$
        \For {$i=L-1, \ldots, 0$}
        \State $\displaystyle \bu^f_{\tau_{i-1}} = \frac{\sqrt{\bar{\alpha}_{\tau_{i-1}}}}{\sqrt{\bar{\alpha}_{\tau_{i}}}} \bu^f_{\tau_{i}} + (\frac{\sqrt{1-\bar{\alpha}_{\tau_{i}}} \sqrt{\bar{\alpha}_{\tau_{i-1}}}}{\sqrt{\bar{\alpha}_{\tau_{i}}}} 
        - \sqrt{1 - \bar{\alpha}_{\tau_{i-1}}})\bs_\theta\left(\bu^f_{\tau_{i}}, \bu^c, \ba^f, {\tau_{i}} \right)$
        \EndFor 
        \For{j=$1, \cdots, t_s$}
        \State Refine $\bu^f_0$ by Gaussian Newton Algorithm~\ref{alg:gn}
        \EndFor \\
        \Return $\bu^f_0$
	\end{algorithmic} 
\end{algorithm}

\section{Numerical experiments}\label{sec:numerical}
In this section, we demonstrate the accuracy and efficiency of PGDM by applying it for solution downscaling in three types of nonlinear PDEs: the nonlinear Poisson equation in both 2D and 3D, the 2D Allen-Cahn equation, and the 2D Navier-Stokes equation. We specifically compare the performance of PGDM against several baseline downscaling techniques, including Fourier Neural Operator (FNO), Cubic Spline Interpolation (CSI), and their enhanced versions that undergo the same number of Gauss-Newton steps as PGDM. Our numerical results indicate that PGDM surpasses all baseline methods in performance, achieving comparable accuracy to high-fidelity solvers while significantly reducing computational costs by more than tenfold. Below, we present in detail the data generation process, the neural network architectures employed, and the numerical results for each test case. 



\subsection{Data generation}

Below we outline the process of generating training and validation data. Let us start with describing the process of  generating the low-fidelity solution, high-fidelity solution, and reference solution. Our investigation covers both  stationary PDEs and time evolutional PDEs. For stationary PDEs, we restrict our attention on the homogeneous Dirichlet boundary condition: 
\begin{equation}\label{eqn:stationary}
    \begin{cases}
    \mL u = a \text{ on } \Omega \\
    u = 0 \text{ on } \partial\Omega.
\end{cases}
\end{equation}
 We employ the finite difference method to discretize the nonlinear differential operator $\mL$, adhering to the specified boundary condition $g(x)$. This discretization, denoted as $\mL_d$, transforms the problem into a nonlinear optimization problem:
\begin{equation}\label{eqn:sta_min}
    \bu \in \arg \min_{\bu} \| \mL_d \bu - \ba \|_2,
\end{equation}
While many nonlinear optimization solvers could potentially be used, in this paper, we opt for the Levenberg–Marquardt (LM) algorithm due to its adaptivity. Specifically, the LM algorithm interpolates between the gradient descent method and the Gauss-Newton method. Throughout the iteration process, it adjusts its behavior, resembling the gradient descent method when the iterates are distant from a local minimum and resembling the Gauss-Newton method when they approach a local minimum. See more details on LM algorithm in Algorithm~\ref{alg:lm} of Appendix \ref{sec:LM}. To generate low-fidelity solution, high-fidelity solution, and reference solution, we execute the the LM algorithm for until the $L^2$-misfit \eqref{eqn:sta_min} decreases below a pre-defined error precision $\eta=$1e-5.


We also consider evolutional PDEs modeled by 
$$
\partial_t u = \mL u,
$$ 
with either Dirichlet or periodic boundary condition.  To ensure the stability while maintaining a reasonable time step size, we adopt an implicit-Euler scheme for evolutional PDEs. More concretely, let $K_t$ be the total number of iteration steps and for $n=0,\cdots, K_t-1$, the approximation solutions $\bu^{n+1}$ are obtained by solving the following optimization problem:
\begin{equation}\label{eqn:min_t}
    \bu^{n+1} \in \arg \min_{\bu} \| (\mI - \Delta t \mL_d) \bu - \bu^n \|_2,
\end{equation}
where $\Delta t$ be the time step size. 
Similar to the static case, we employ the LM algorithm as our numerical solver for $\eqref{eqn:min_t}$ with the stopping criterion set to be that  the $L^2$-misfit is below $\eta=$1e-5. For the two evolutional PDEs considered in the paper, namely the Allen-Cahn and Navier-Stokes equations, we consider spatial super-resolution only.  Specifically, we employ spatial mesh grids denoted as $K_c$ and $K_f$, along with a time step size of $\Delta t=0.05$ for both the low-fidelity and high-fidelity solutions. The reference solutions for the Allen-Cahn equation are computed on a spatial mesh size of $2K_f$ and a time step size of $\Delta t = 0.025$. For the reference solutions of Navier-Stokes equation, we adopt a Crank-Nicolson scheme commonly used in literature, such as the one utilized in \cite{li2020fourier}.  Specifically, we set the time step size to be $\Delta t =$5e-5 and utilize a spatial mesh grid size of $2K_f$. It is important to note that there is a trade-off between using an implicit scheme with a larger time step and using a semi-implicit scheme with a smaller time step. The former allows for a larger time step, leading to faster evolution, but it introduces an error of $O(\Delta t)$.

In the data preparation step, we generate the source terms and the initial conditions from the Gaussian random field $\mN(0, (-\Delta + b^2 \mI)^{-c})$, where $b$ and $c$ are two hyperparameters adjusting the length-scale and smoothness of the field. For comparison purposes, all comparisons in this section are conducted at the resolution of the high-fidelity solution, indicated by a spatial mesh size $K_f$. Let us introduce three classical solvers for solving the nonlinear systems \eqref{eqn:sta_min} or \eqref{eqn:min_t}. The coarse solver generates low-fidelity solutions by solving these systems on a coarser mesh grid $K_c$ using LM algorithm and subsequently enhances resolution through cubic spline interpolation. The fine solver produce high-fidelity solutions by directly solves the nonlinear systems employing a finer spatial mesh grid $K_f$. Additionally, the reference solver utilize even finer spatial mesh grid $2K_f$ coupled with a significantly finer time step size, followed by downsampling to match the resolution of fine solver. In each of the following numerical examples, we employ the aforementioned methodology to generate $N$ training sample tuples $\{ \bu_i^c, \bu_i^f, \ba_i \}_{i=1}^N$ and $M$ testing sample tuples $\{ \bu_j^c, \bu_j^f, \bu^r_j, \ba_j \}_{j=1}^M$. Here  $\ba$ is generated by sampling from the Gaussian random field followed by restriction on the grids, $\bu^c$ is obtained by the CSI solver,  $\bu^f$ is obtained by the fine solver and  $\bu^r$ is obtained by the reference solver. Given our focus on scenarios with limited data, we set $N$ in this paper to be as small as $N=30$. For a better illustration, we summarize the previously mentioned notations and hyperparameters along with their definitions in Table~\ref{tab:notation} in Appendix. Additionally, the detailed descriptions of the neural network architecture and hyperparameters used for the diffusion models and FNO can be found in Appendix~\ref{sec:archi}.

\subsection{Nonlinear Poisson equation}
\begin{figure}[h!]
    \centering
    \includegraphics[width=1\textwidth, height=5cm, keepaspectratio]{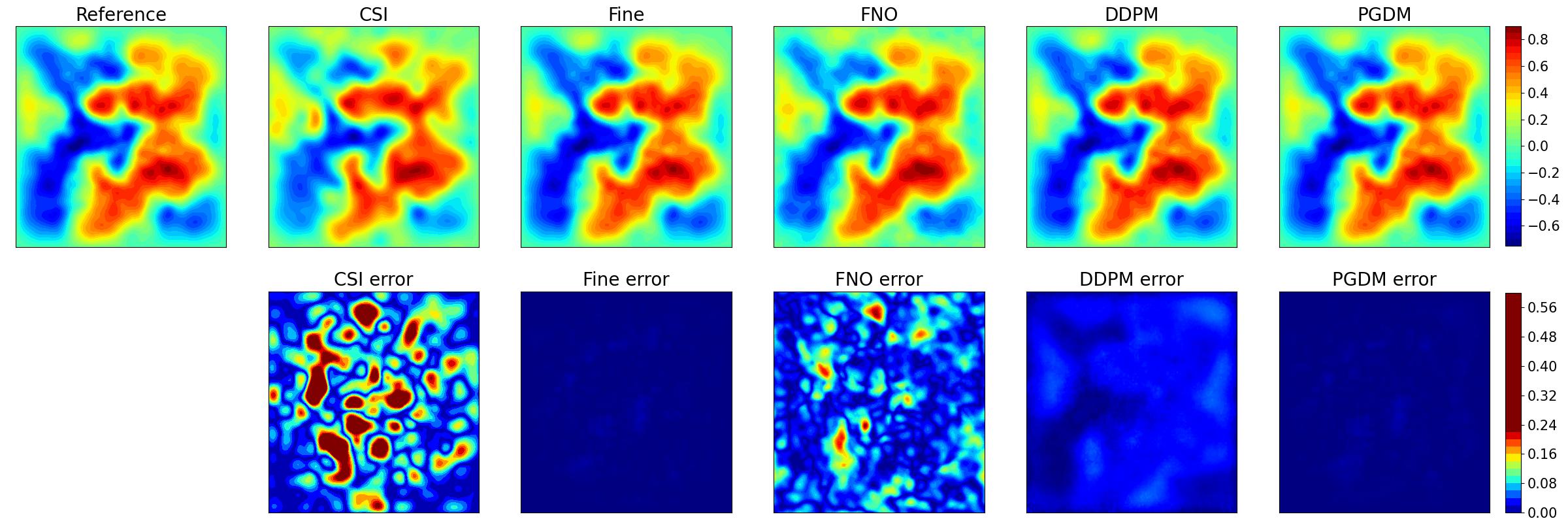}
    \caption{2D nonlinear Poisson: Predictions and absolute errors generated by different solvers with $c=1.6$ and $N=100$ training samples.}
    \label{fig:P_2d_1p6}
\end{figure}

\begin{figure}[h!]
    \centering
    \includegraphics[width=1\textwidth, height=5cm, keepaspectratio]{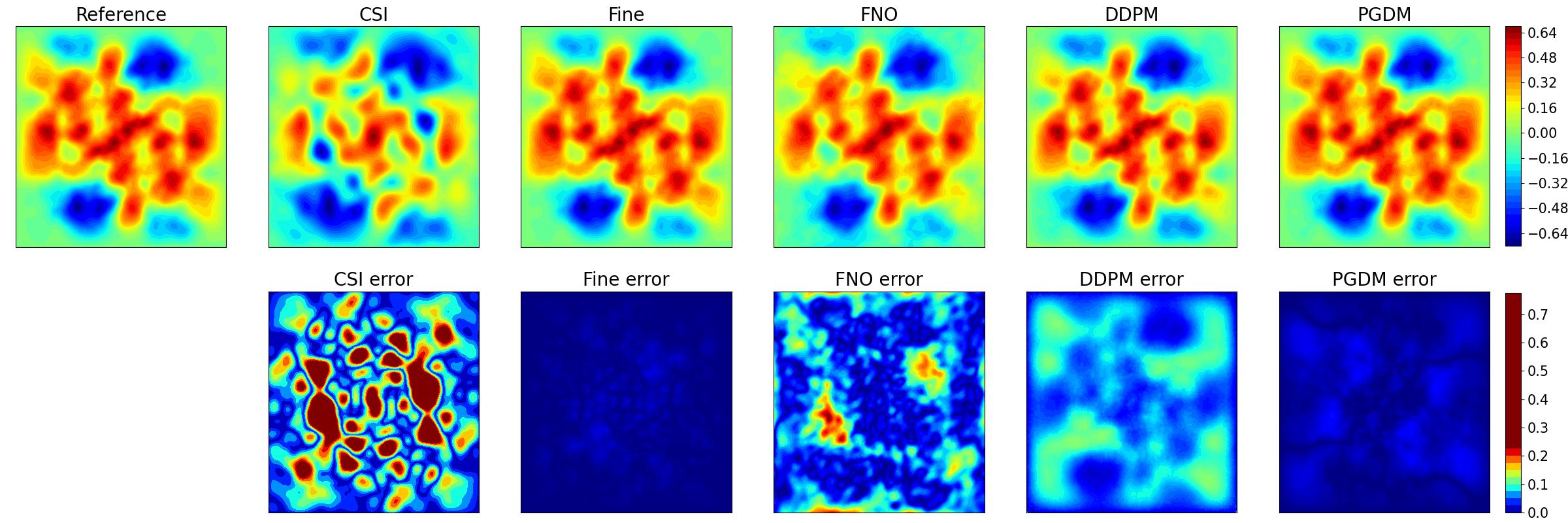}
    \caption{2D nonlinear Poisson: Predictions and corresponding absolute errors generated by different solvers with $c=1.2$ and $N=100$ training samples.}
    \label{fig:P_2d_1p2}
\end{figure}

\begin{figure}[h!]
    \centering
    \includegraphics[width=1\textwidth, height=5cm, keepaspectratio]{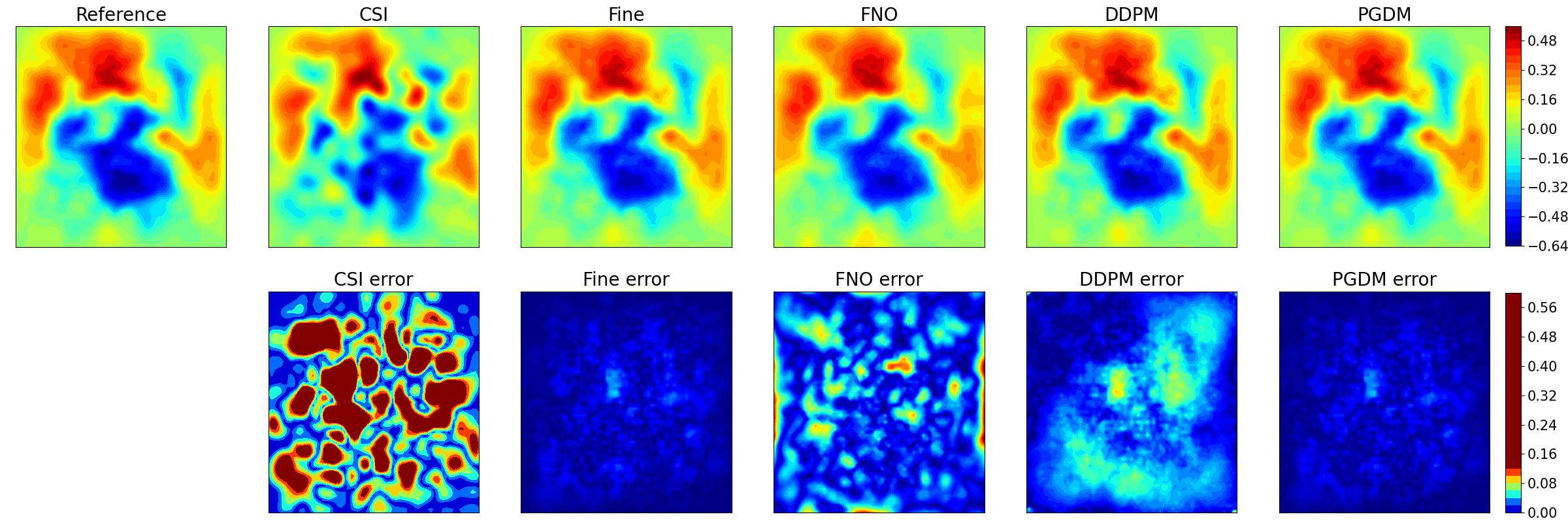}
    \caption{3D nonlinear Poisson: Predictions and corresponding absolute errors generated by different solvers with $c=1.6$ and $N=100$ training
samples.}
    \label{fig:P_3d_1p6}
\end{figure}

\begin{figure}[h!]
    \centering
    \includegraphics[width=1\textwidth, height=5cm, keepaspectratio]{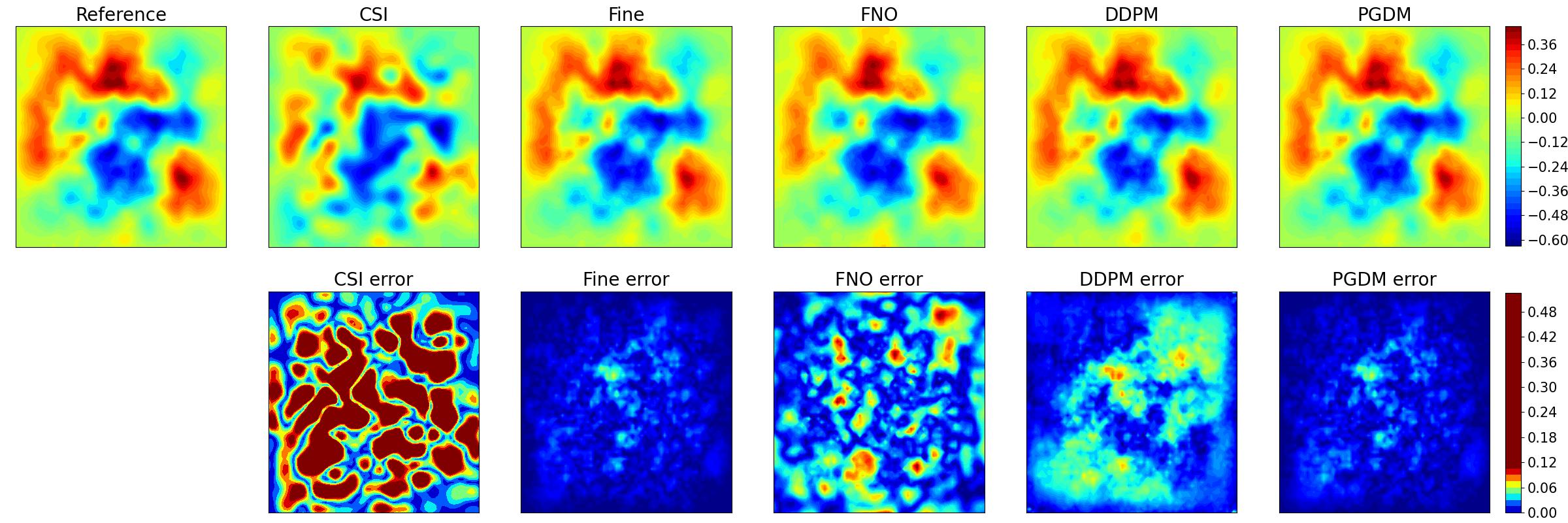}
    \caption{3D nonlinear Poisson: Predictions and corresponding absolute errors generated by different solvers with $c=1.4$ and $N=100$ training
samples.}
    \label{fig:P_3d_1p4}
\end{figure}

The first example is the nonlinear Poisson equation with zero Dirichlet boundary condition,
\begin{equation}\label{eqn:nonlp}
    \begin{split}
     - 0.0005 \Delta u(x) + u(x)^3 &= a(x), \quad x\in{(0,1)}^{d}, \; \\
    u(x) &= 0,  \quad \quad x \in \partial (0, 1)^d.
\end{split}
\end{equation}
Here, $d$ indicates the physics dimensionality and $a(x)$ denotes the source term, which is sampled from a Gaussian random fields described by $\mN(0, (-\Delta + 49 \mI)^{-c})$, where the inverse Laplacian is equipped with zero boundary condition. Our investigation spans various values of $c$ and the size of training set $N$. For 2D cases, i.e. $d=2$, we select $K_c = 16$, $K_f=128$ as the mesh grid sizes for the coarse solver and the fine solver, respectively. In 3D scenarios, we select $K_c = 16$, $K_f=64$ for the mesh grid size of the coarse and the fine solver respectively. The performance of different solvers under these conditions for the 2D scenario is detailed in Table~\ref{tab:compare_DM_Nonl_2d}. Solutions computed from various solvers in 2D with $c=1.6$ and $N=100$ are depicted in Figure~\ref{fig:P_2d_1p6}, while the corresponding solutions for $c=1.2$ and $N=100$ are illustrated in Figure~\ref{fig:P_2d_1p2}. Additionally, the performance of different solvers under analogous conditions in 3D settings is comprehensively detailed in Table~\ref{tab:compare_DM_Nonl_3D}. Visual comparisons for 3D solver outputs corresponding to $c=1.2$ and $N=100$ are also presented are presented in Figure~\ref{fig:P_3d_1p6}, and the results for $c=1.4$ and $N=100$ are shown in Figure~\ref{fig:P_3d_1p4}. As demonstrated in Table~\ref{tab:compare_DM_Nonl_2d} and Table~\ref{tab:compare_DM_Nonl_3D}, PGDM maintains the same level of accuracy as the fine solver while reducing computational time by a significant factor of ten.

\begin{table}[h!]
\begin{center}
\begin{tabular}{|c|c|c|c|c||c|}
\hline
 & \begin{tabular}{@{}c@{}}$N=30$\\ $c=1.6$ \end{tabular}    & \begin{tabular}{@{}c@{}}$N=100$\\ $c=1.6$ \end{tabular}  & \begin{tabular}{@{}c@{}}$N=30$\\ $c=1.2$ \end{tabular}  & \begin{tabular}{@{}c@{}}$N=100$\\ $c=1.2$ \end{tabular}  & Time  \\
\hline
\hline
 CSI  & 2.97e-1  & 2.97e-1  & 5.82e-1 & 5.82e-1  & 3.05e-1 \\
 \hline
 Fine  & 3.69e-3 & 3.69e-3 & 1.36e-2 & 1.36e-2 &  8.73e0 \\
 \hline
 \hline
  FNO & 2.36e-1 & 1.73e-1 & 3.36e-1 & 2.22e-1  & 1.66e-1  \\
 \hline
 DDPM & 8.74e-2 & 6.48e-2 & 1.44e-1 & 1.36e-1 &  3.22e0 \\
  \hline
 DDIM & 1.18e-1 & 6.83e-2 & 1.66e-1 & 1.38e-1 &  6.15e-1 \\
  \hline
  \hline
Coarse+GN  & 6.17e-2 & 6.17e-2 & 1.91e-1 & 1.91e-1  & 6.27e-1  \\
\hline
FNO+GN & 3.93e-2 & 2.28e-2 & 1.07e-1 & 4.67e-2 & 6.14e-1  \\
\hline
 PGDM  & \textbf{1.31e-2} & \textbf{5.20e-3} & \textbf{2.64e-2} & \textbf{2.01e-2} & \textbf{1.29e0}  \\
 \hline
\end{tabular}
\end{center}
\caption{Comparison of relative $L^2$-error for 2D nonlinear Poisson equation at 8x super-resolution on $M=30$ testing examples. The last column shows the average computational time over $M=30$ realizations of different solutions.}
\label{tab:compare_DM_Nonl_2d}
\end{table}


\begin{table}[h!]
\begin{center}
\begin{tabular}{|c|c|c|c|c||c|}
\hline
 & \begin{tabular}{@{}c@{}}$N=30$\\ $c=1.6$ \end{tabular}    & \begin{tabular}{@{}c@{}}$N=100$\\ $c=1.6$ \end{tabular}  & \begin{tabular}{@{}c@{}}$N=30$\\ $c=1.4$ \end{tabular}  & \begin{tabular}{@{}c@{}}$N=100$\\ $c=1.4$ \end{tabular}  & Time  \\
\hline
\hline
 CSI  & 4.69e-1  & 4.69e-1& 6.72e-1 & 6.72e-1 & 2.18e-1  \\
 \hline
 Fine  & 2.66e-2 & 2.66e-2 & 4.74e-2 & 4.74e-2 & 1.60e2 \\
 \hline
 \hline
  FNO & 2.29e-1 & 1.71e-1 & 2.73e-1 & 2.03e-1 & 6.63e-1 \\
 \hline
 DDPM & 1.13e-1 & 1.10e-1 & 1.12e-1& 1.10e-1 & 3.23e1 \\
  \hline
 DDIM & 1.46e-1 & 1.33e-1 & 1.43e-1 & 1.36e-1 & 6.68e0 \\
  \hline
  \hline
Coarse+GN  & 9.61e-2 & 9.61e-2  & 1.61e-1 & 1.61e-1 & 9.91e0\\
\hline
FNO+GN & 3.97e-2 & 2.99e-2 & 6.26e-1 & 5.19e-2 & 1.03e1 \\
\hline
 PGDM  & \textbf{2.79e-2} & \textbf{2.77e-2} &\textbf{4.77e-2} & \textbf{4.74e-2} & \textbf{1.63e1} \\
 \hline
\end{tabular}
\end{center}
\caption{Comparison of relative $L^2$ error for 3D nonlinear Poisson equation at 4x super-resolution on $M=30$ testing set. The last column documents the average computational time over $M=30$ realizations of different solutions.}
\label{tab:compare_DM_Nonl_3D}
\end{table}

\newpage

\subsection{2D Allen-Cahn equaiton}
\begin{figure}[h!]
    \centering
    \includegraphics[width=0.9\textwidth, height=5cm, keepaspectratio]{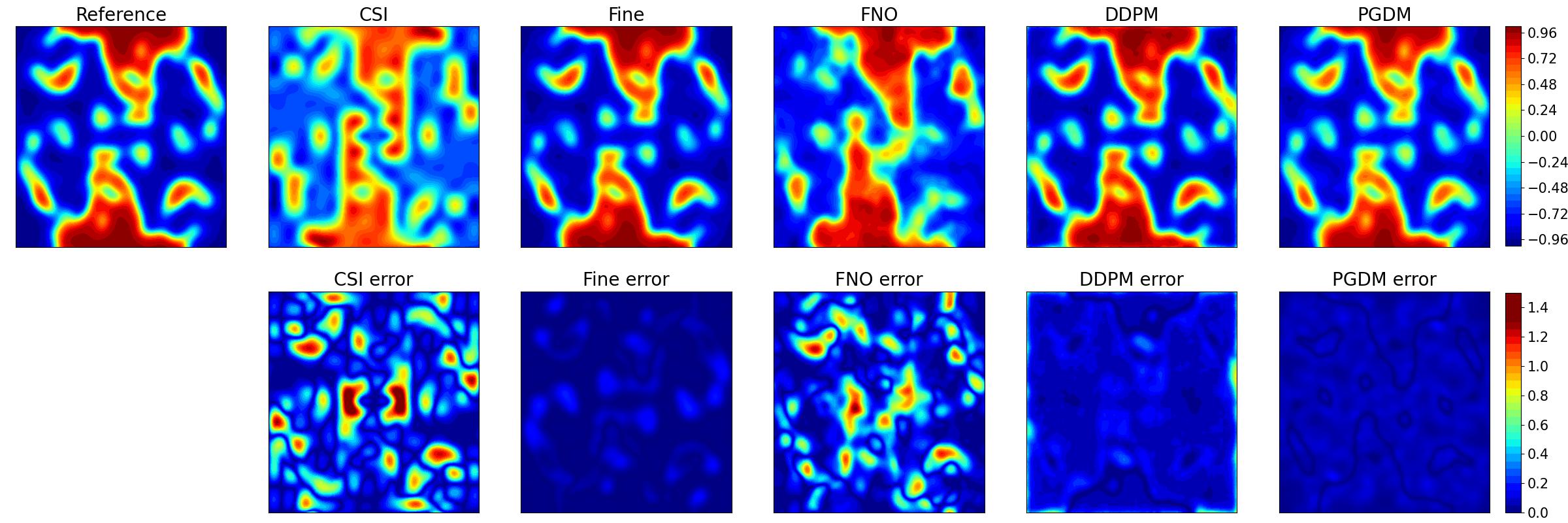}
    \caption{2D Allen-Cahn: Predictions at $t=0.5$ and corresponding absolute errors generated by different solvers with $\gamma=5$, $c=1.6$ and $N=30$.}
    \label{fig:P_ac_1p6}
\end{figure}

\begin{figure}[h!]
    \centering
    \includegraphics[width=1\textwidth, height=5cm, keepaspectratio]{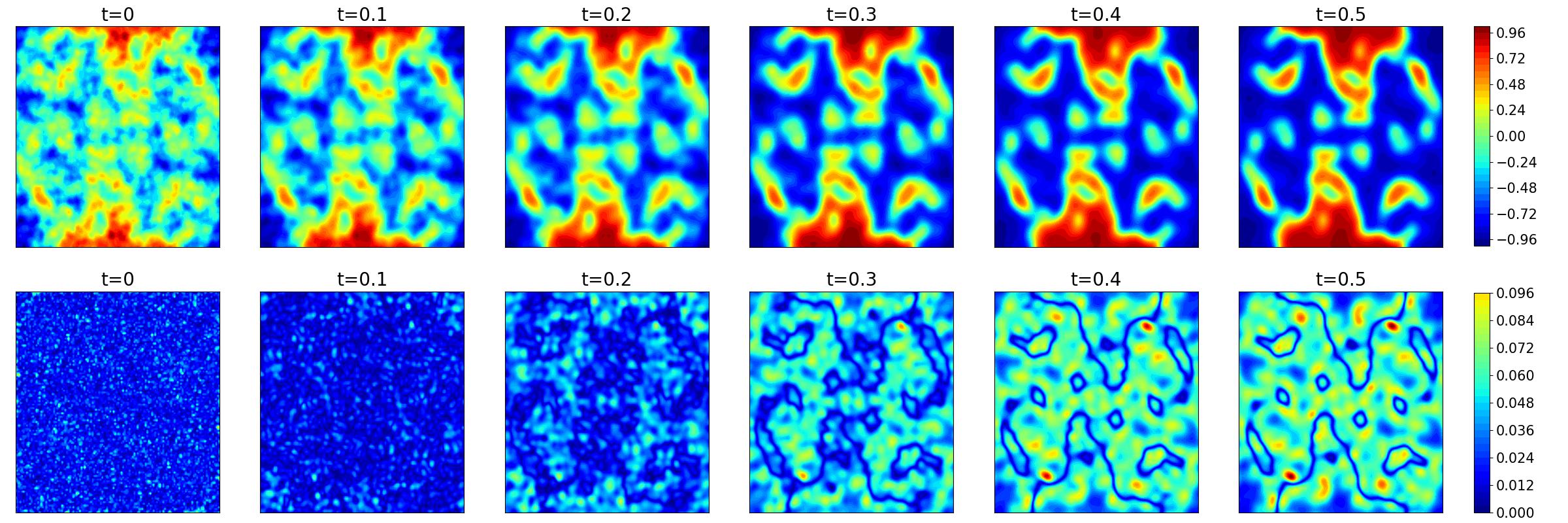}
    \caption{2D Allen-Cahn: Snapshots of evolution of PGDM  and the corresponding absolute errors compared to the reference solution with $\gamma=5$, $c=1.6$ and $N=30$.}
    \label{fig:P_ac_1p6_evo}
\end{figure}

Consider 2D Allen-Cahn equaiton with periodic boundary condition:
\begin{equation}\label{eqn:ac}
    \begin{split}
    \partial_t u(x) &= \kappa \Delta u(x) + \gamma u(x) \Big(\frac{1}{4}-u(x)^2\Big), ~ ~ x \in (0,1)^2, t \in (0, 0.5], \; \\
    u(0, x) &= u_0(x), ~ ~ x \in (0,1)^2.
\end{split}
\end{equation}
Here the diffusion coefficient set to $\kappa=$1e-3 and the reaction coefficient set to $\gamma=5$. The initial conditions $u_0(x)$ draw from Gaussian random field $\mN(0, (-\Delta + 49 \mI)^{-c})$ where the inverse Laplacian $\Delta$ is applied with periodic boundary conditions. We explore different values of reaction coefficient $\gamma$ and various $c$. The time step size is set to $\Delta t = 0.05$, with the total number of steps $K_t$ set to 10. For computational mesh grids, sizes are set at $K_c = 16$ for the coarse solver and $K_f = 128$ for the fine solver, respectively. The performance of various solvers across these settings is detailed in Table~\ref{tab:compare_DM_AC}. Predictions at $t=0.5$ of different solvers with $c=1.6$, $\gamma=5$ and $N=30$ are presented in Figure~\ref{fig:P_ac_1p6}, and the snapshots of predictions of PGDM are shown in Figure~\ref{fig:P_ac_1p6_evo}.


\begin{table}[h!]
\begin{center}
\begin{tabular}{|c|c|c||c||c|c||c||}
\hline
 & \begin{tabular}{@{}c@{}}$\gamma=1$\\ $c=1.6$ \end{tabular}    & \begin{tabular}{@{}c@{}}$\gamma=1$\\ $c=1.2$ \end{tabular}  & Time & \begin{tabular}{@{}c@{}}$\gamma=5$\\ $c=1.6$ \end{tabular}  & \begin{tabular}{@{}c@{}}$\gamma=5$\\ $c=1.2$ \end{tabular}  & Time  \\
\hline
\hline
 CSI & 4.16e-1 & 7.84e-1 & 5.59e-2  & 4.57e-1   & 8.37e-1 &  2.64e-1 \\
 \hline
 Fine & 2.67e-2 & 6.84e-2 & 1.71e1 & 4.73e-2  & 7.73e-2  & 5.84e1 \\
 \hline
 \hline
  FNO & 4.57e-1 & 7.27e-1 & 8.21e-2 & 4.61e-1   & 7.84e-1 & 8.21e-2 \\
 \hline
 DDPM & 1.66e-1 & 1.71e-1 & 1.42e1 & 2.05e-1   & 1.89e-2 & 1.42e1 \\
  \hline
 DDIM & 1.82e-1 & 3.2e-1 & 2.84e0 & 2.11e-1   & 3.54e-1 & 2.84e0 \\
  \hline
  \hline
Coarse+GN & 8.73e-2 & 1.91e-1 & 3.18e-1 & 1.06e-1   & 1.94e-1 & 4.87e-1 \\
\hline
FNO+GN & 1.03e-1 & 1.87e-1 & 3.42e-1 & 1.29e-1   & 2.05e-1 & 3.42e-1 \\
\hline
 PGDM & \textbf{5.96e-2} & \textbf{9.09e-2} & \textbf{3.11e0} & \textbf{6.67e-2}   & \textbf{1.05e-1} & \textbf{4.23e0} \\
 \hline
\end{tabular}
\end{center}
\caption{Comparison of relative $L^2$ error for 2D Allen-Cahn equation at 8x super-resolution on $N=30$ training set and $M=20$ testing set. The fourth and the last column show the averaged computational time over $M=20$ realizations of different solutions.}
\label{tab:compare_DM_AC}
\end{table}

\newpage

\subsection{2D Navier-Stokes equation}
\begin{figure}[h!]
    \centering

    \includegraphics[width=1\textwidth, height=5cm, keepaspectratio]{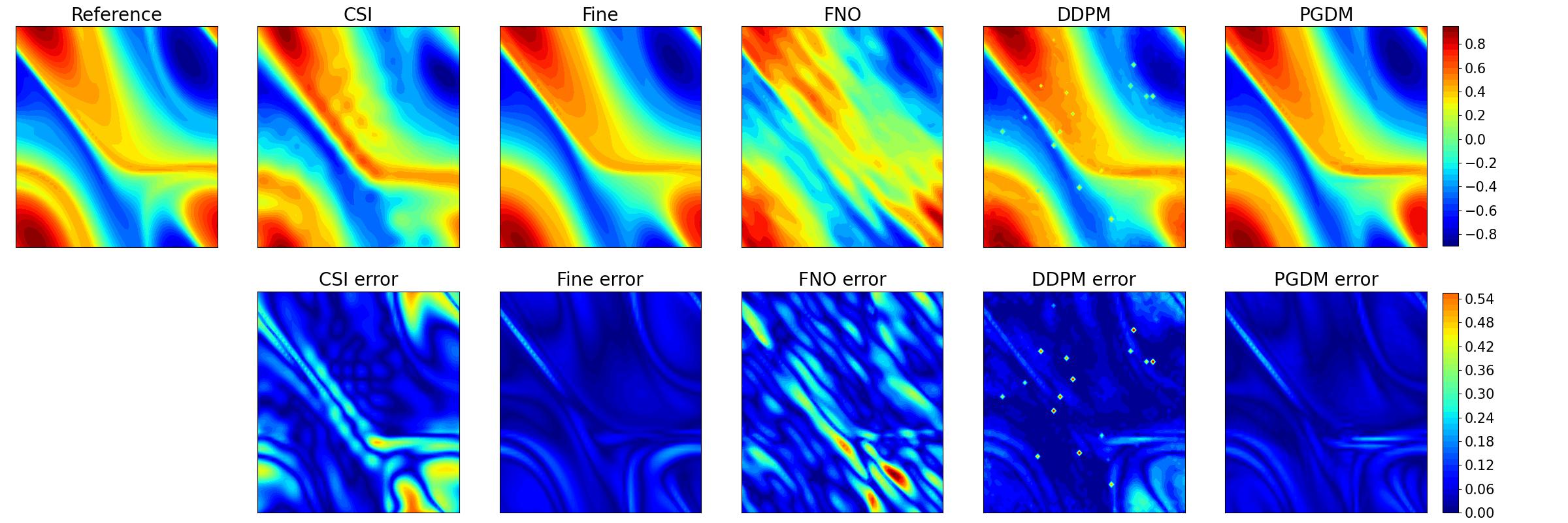}

    \caption{2D Navier-Stokes: Predictions at $t=2$ and the corresponding absolute errors generated by different solvers with $\nu=$2e-4 and $N=30$.}
    \label{fig:ns_nu_2en4}
\end{figure}

\begin{figure}[h!]
    \centering

    \includegraphics[width=1\textwidth, height=5cm, keepaspectratio]{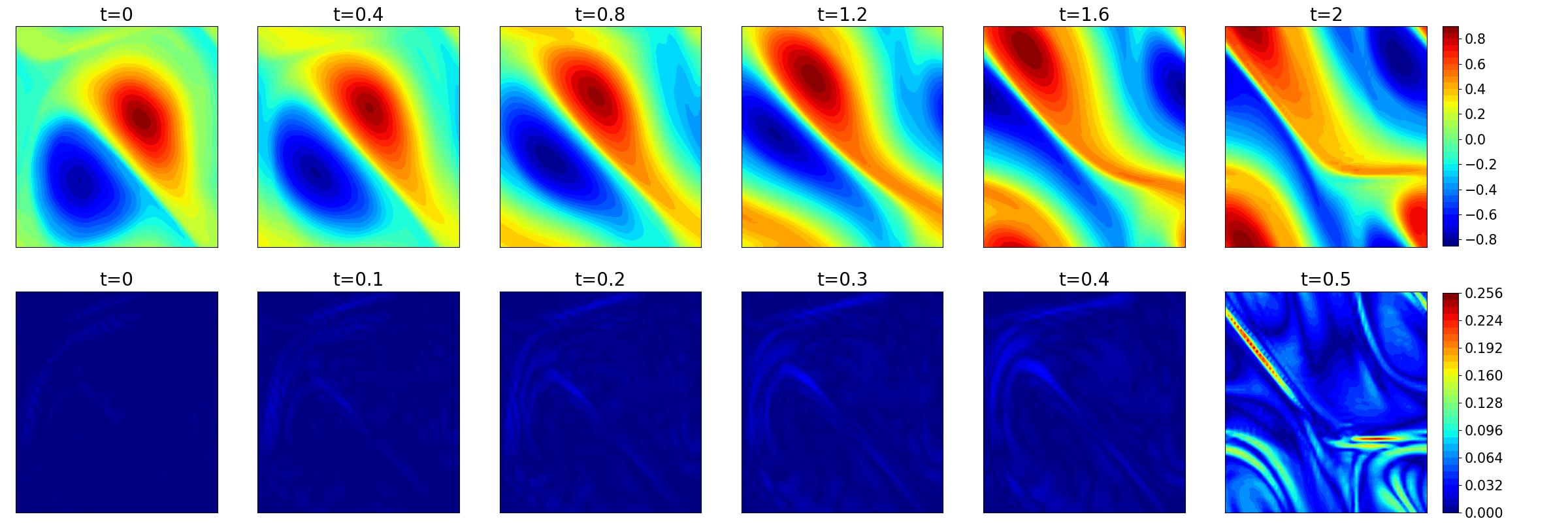}
    
    \caption{2D Navier-Stokes: Snapshots of evolution of PGDM and the corresponding absolute errors to the reference solution with $\nu=$2e-4 and $N=30$.}
    \label{fig:ns_nu_2en4_2}
\end{figure}

\begin{figure}[h!]
    \centering
    \includegraphics[width=1\textwidth, height=5cm, keepaspectratio]{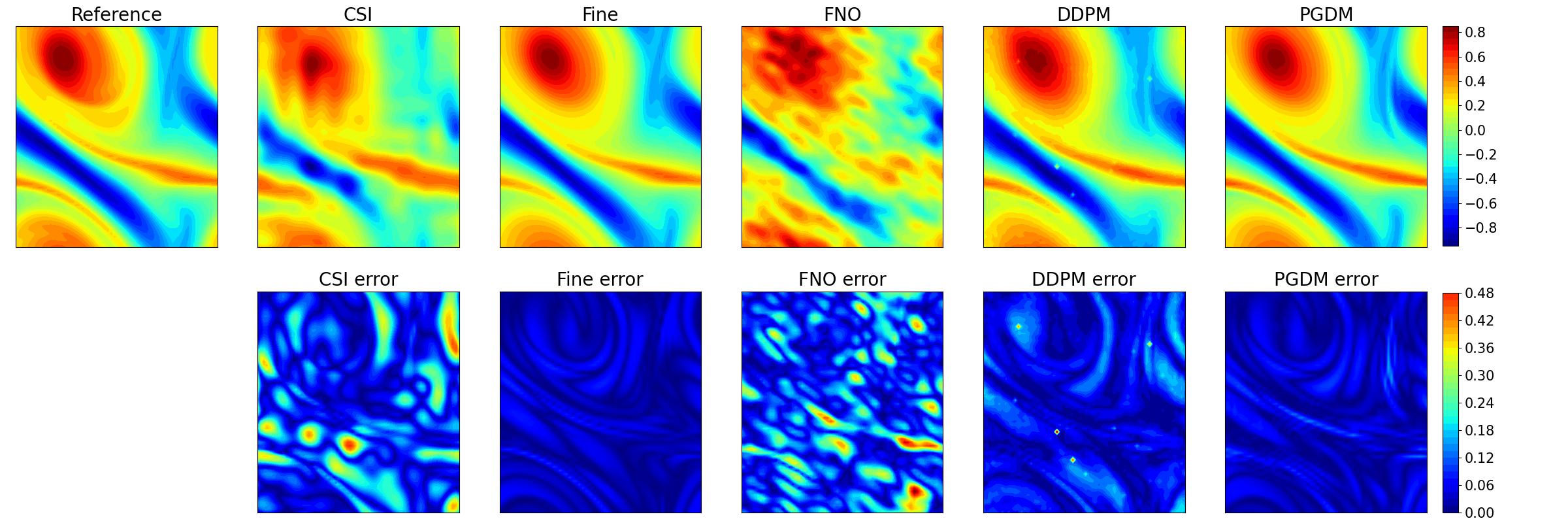}
    \caption{2D Navier-Stokes: Predictions at $t=2$ and the corresponding absolute errors generated by different solvers with $\nu=$1e-4 and $N=30$.}
    \label{fig:ns_nu_1en4}
\end{figure}

\begin{figure}[h!]
    \centering

    \includegraphics[width=1\textwidth, height=5cm, keepaspectratio]{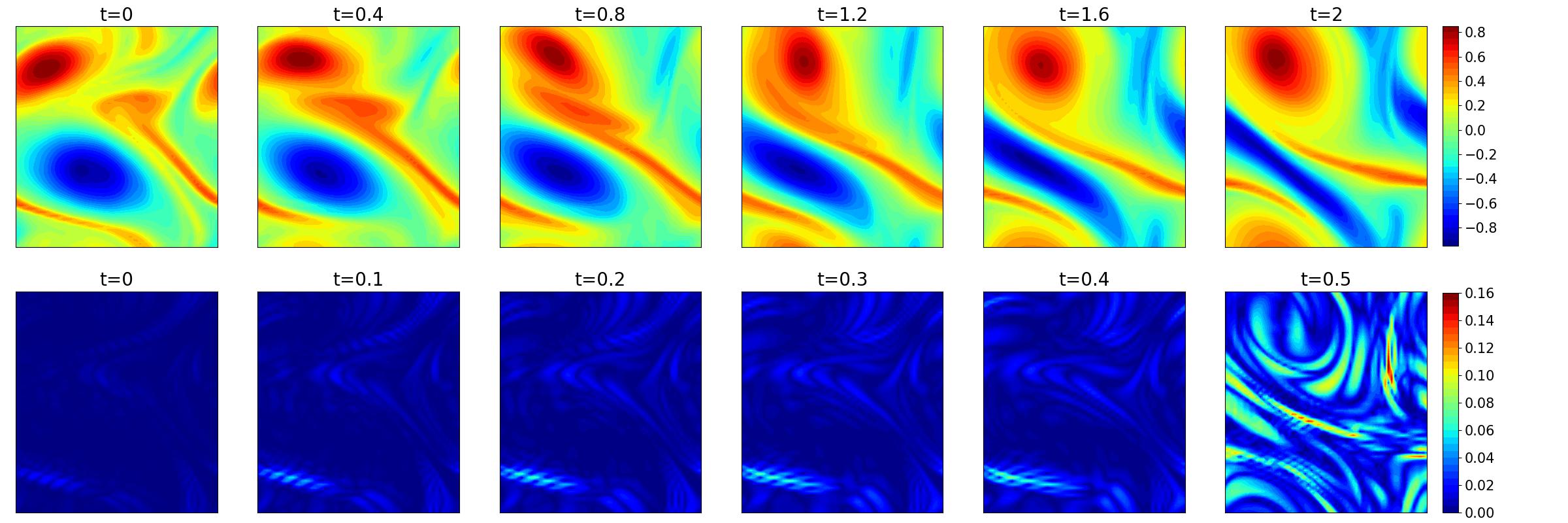}
    
    \caption{2D Navier-Stokes: Snapshots of evolution of PGDM  and the corresponding absolute errors to the reference solution with $\nu=$1e-4 and $N=30$.}
    \label{fig:ns_nu_1en4_2}
\end{figure}

Consider the 2D Navier-Stokes equation in the vorticity form with periodic boundary condition:
\begin{equation} \label{eqn:ns}
\begin{split}
\partial_{t} w(t, x)+ \mu u(t, x) \cdot \nabla w(t, x) &=\nu \Delta w(t, x)+  f(x), \quad x \in (0,1)^2, t \in (0, 2], \\
w(0, x) &= w_0(x), \qquad  \qquad \qquad x \in (0,1)^2.
\end{split}
\end{equation}
The transportation coefficient is set to $\mu=4$, and the forcing term is selected as $f(x) = 0.2 (\sin(2 \pi (x + y)) + \cos(2 \pi (x + y)))$. To generate initial condition, we draw functions from the same Gaussian random field $\mN(0, (-\Delta + 25 \mI)^{-5})$ where the inverse Laplacian $\Delta$ is applied with periodic boundary conditions. These functions are subsequently evolved using the reference solver for two seconds. The time step size is set to $\Delta t = 0.05$, with the total number of steps set to $K_t = 40$. The mesh grid sizes are set to $K_c = 16$ for the coarse solver and $K_f = 64$ for the fine solver, respectively. We fix the training set size $N$ at 30 and examine three different viscosity coefficients $\nu =$ 2e-4, 1e-4, 2e-5. The performance of various solvers across these settings is detailed in Table~\ref{tab:compare_DM_NS}. Predictions at $t=2$ of different solvers with $\nu=$2e-4 are presented in Figure~\ref{fig:ns_nu_2en4}, and the snapshots of predictions of PGDM are shown in Figure~\ref{fig:ns_nu_2en4_2}; Predictions at $t=2$ of different solvers $\nu=$1e-4 are presented in Figure~\ref{fig:ns_nu_1en4}, and the snapshots of predictions of PGDM are shown in Figure~\ref{fig:ns_nu_1en4_2};

\begin{table}[h!]
\begin{center}
\begin{tabular}{|c|c|c|c||c|}
\hline
 &  $\nu =$ 2e-4   &  $\nu =$  1e-4  &  $\nu =$ 2e-5  & Time  \\
\hline
\hline
 CSI & 1.54e-1 & 1.63e-1 & 2.35e-1 & 1.38e0 \\
 \hline
 Fine & 4.23e-2 & 6.51e-2 & 1.19e-1 & 2.26e2 \\
 \hline
 \hline
  FNO & 3.43e-1 & 3.65e-1 & 4.88e-1 & 1.38e-2 \\
 \hline
 DDPM & 8.09e-2 & 1.06e-1 & 1.56e-1 & 8.99e-1\\
  \hline
 DDIM & 9.59e-2 & 1.43e-1 & 2.01e-1 & 2.32e-1\\
  \hline
  \hline
Coarse+GN & 5.42e-2 & 7.98e-2 & 1.33e-1 & 6.28e0 \\
\hline
FNO+GN & 1.22e-1 & 1.41e-1 & 1.84e-1 & 6.13e0 \\
\hline
 PGDM & \textbf{4.23e-2} & \textbf{6.51e-2} & \textbf{1.19e-1}  & \textbf{6.34e0} \\
 \hline
\end{tabular}
\end{center}
\caption{Comparison of relative $L^2$ error for 2D Navier-Stokes equation at 4x super-resolution on $M=20$ testing set. The last column documents the average computational time over $M=20$ realizations of different solutions.}
\label{tab:compare_DM_NS}
\end{table}

\section{Conclusion}
We propose a data-driven surrogate method called PGDM for accelerating the computation of (nonlinear) PDEs. PGDM first generates a high-fidelity solution conditional on a low-resolution input, followed by a mild refinement of the former with a  PDE solver on the fine-grid.   Our numerical results show that  PGDM can produce high-fidelity solutions that are comparable to those generated by fine-scale solvers, while requiring very limited training data (as few as 30 instances). More importantly, we demonstrate that our PGDM also significantly reduce the computational time, especially in 3D examples, where we observe a tenfold decrease compared to fine-scale solvers.

\section*{Acknowledgement}
YL thanks the support from the National Science Foundation through the award DMS-2343135 and the support from the Data Science Initiative at University of Minnesota through a MnDRIVE DSI Seed Grant.

\appendix
\section*{Appendix}

\begin{table}
\begin{center}
\begin{tabular}{c c}
\hline
 Notation & Meaning \\
 \hline
 $K_c$ & Uniform grid size of coarse solver \\
 $\bu^c$ & CSI solution \\
 $K_f$ & Uniform grid size of fine solver \\
 $\bu^f$ & Fine solution \\
 $K_t$ & Number of evolution steps \\
 $\Delta t$ & Time step size in evolution problem \\
 $T$ & Total time steps in DDPM \\
 $\tau$ & Sequence of skipped time steps in DDIM \\
 $\{ \beta_i \}_{i=0}^T$ & Scale of Noise in DDPM \\
 $\{ \alpha_i \}_{i=0}^T$ & Hyperparameter in DDPM \\
 $N$ & Number of training samples \\
 $M$ & Number of testing samples \\
 $b,c$ & Hyperparameters of Gaussian random fields \\
 \hline
\end{tabular}
\end{center}
\caption{Table of notations}
\label{tab:notation}
\end{table}

\section{Neural Networks Architecture and Hyperparameters}\label{sec:archi}
Our diffusion models are based on the DDPM architecture \cite{ho2020denoising}, which uses U-Net \cite{ronneberger2015u} as the backbone. During our experiments, we omit the use of self-attention, resulting in significant reductions in training time while maintaining similar sample quality. The base channel count, the list of Down/Up channel multipliers and the list of middle channel refer to the hyperparameters of the U-Net, which is detailed in Table~\ref{tab:DM_model_parameters}. To accelerate sampling process using DDIM, we take skipped time steps  $\tau$ be $[1, 5, 10, 15, 20, 25, \cdots, T-5, T]$. The linear noise schedule is configured from $\beta_0 = 0.0001$ to $\beta_{T}=0.02$. During training, we utilize the Adam optimizer with a dynamic learning rate that linearly decays every 5000 steps with a decay rate of 0.05. The total number of training steps is set to 10000.

\begin{table}
\begin{center}
\begin{tabular}{|c|c|c|c|c|}
\hline
 & \begin{tabular}{@{}c@{}} 2D Nonlinear\\ Poisson \end{tabular} & \begin{tabular}{@{}c@{}} 3D Nonlinear\\ Poisson \end{tabular} & \begin{tabular}{@{}c@{}} 1Dt + 2Dx \\ Allen-Cahn \end{tabular} & \begin{tabular}{@{}c@{}} 1Dt + 2Dx \\ Navier-Stokes \end{tabular} \\
  \hline
  Target resolution  & $128 \times 128$ & $64 \times 64 \times 64$ & $10 \times 128 \times 128$ & $40 \times 64 \times 64$ \\
\hline
 \begin{tabular}{@{}c@{}} Super-resolution \\ scale \end{tabular} & 8  & 4 & 8 & 4 \\
\hline
\hline
 Timesteps $T$ & 400  & 400 & 400 & 200 \\
\hline
 Base Channels & 128  & 128 & 256 & 128 \\
 \hline
  \begin{tabular}{@{}c@{}} Down/Up \\ Channel multipliers \end{tabular} & 1,2,4,8  & 1,2,4,8 & 1,2,4,8 & 1,2,4,8 \\
 \hline
  \begin{tabular}{@{}c@{}} Middle \\ Channel list \end{tabular} & [512, 512]  & [1024, 1024] & [1024, 1024] & [1024, 1024] \\
 \hline
 \hline
 Batch size & 4  & 2 & 2 & 2 \\
 \hline
\end{tabular}
\end{center}
\caption{Table of DDPM hyperparameters}
\label{tab:DM_model_parameters}
\end{table}

The architecture of FNO follows that described in \cite{li2020fourier}. The number of lifting channels, number of FFT truncation modes, and number of Fourier layers for different examples are specified in Table~\ref{tab:FNO_model_parameters}. During training, we utilize the Adam optimizer with a dynamic learning rate that linearly decays every 5000 steps with a decay rate of $0.05$. Training continues until the loss drops below 1e-6 or reaches the maximum iteration number of 50000.

Our model training were performed on an NVIDIA RTX 3070 graphics card, while predictions and refinements with Gaussian-Newton were executed on an AMD Ryzen 7 3700X processor.

\begin{table}
\begin{center}
\begin{tabular}{|c|c|c|c|}
\hline
 & \text{layers} & \text{modes} & \text{lifting channel}   \\
\hline
 Nonlinear Poisson 2D &2 & 16 & 32 \\
 \hline
  Nonlinear Poisson 3D & 4 & 12 & 64  \\
 \hline
   Allen-Cahn $1Dt+2Dx$ & 4 & 16 & 64  \\
 \hline
Navier-stokes $1Dt+2Dx$ & 4 & 16 & 64  \\
 \hline
\end{tabular}
\end{center}
\caption{Table of FNO hyperparameters}
\label{tab:FNO_model_parameters}
\end{table}

\section{Levenberg–Marquardt algorithm} \label{sec:LM}
In this part, we present the Levenberg–Marquardt (LM) algorithm for solving the nonlinear optimization problem \eqref{eqn:sta_min} and \eqref{eqn:min_t} in Algorithm~\ref{alg:lm}. In all of our numerical experiments, we fix $\lambda=0.5$ and $\eta=$1e-5.
\begin{algorithm}[h]
	\caption{Levenberg–Marquardt algorithm} 
    \label{alg:lm}
	\begin{algorithmic}[1]
        \Require Initial guess $\bu_0$, the source term $\ba$, the discretization of operator $\mL_d$, initial damping  parameter $\lambda$ and stopping criterion $\eta$.
        \State Let $\bu=\bu_0$
        \Repeat
        \State Compute residual $\boldsymbol{r} = {\mL_d} \bu - \ba$
        \State Compute Jacobian matrix $J = \frac{\partial {\mL_d} \bu}{\partial \bu}$
        \State Solve linear system $[J^TJ + \lambda \text{diag}(J^TJ)] \boldsymbol{\delta} = J^T \boldsymbol{r}$ for $\boldsymbol{\delta}$
        \State Update $\bu_{new} = \bu + \delta$
        \If{$\| {\mL_d} \bu_{new} - \ba \| > \| {\mL_d} \bu - \ba \|$} 
        \State $\lambda = \lambda * 2$ \Else \State $\lambda = \lambda / 2$ \EndIf
        \State $\bu = \bu_{new}$
        \Until{$\| \boldsymbol{r} \|< \eta$} \\
    \Return{$\bu$}
	\end{algorithmic} 
\end{algorithm}
\section{Gaussian-Newton Algorithm}
To refine the solution obtained from the coarse solver, diffusion model and the FNO, we introduce the one-step Gaussian-Newton refinement process, outlined in Algorithm~\ref{alg:gn}.
\begin{algorithm}
	\caption{One step Gaussian Newton update} 
    \label{alg:gn}
	\begin{algorithmic}[1]
 \Require The $\bu$ to be refined, the source term $\ba$, and the discretization of operator $\mL_d$.
        \State Compute Jacobian matrix $J = \frac{\partial {\mL_d} \bu}{\partial \bu}$
        \State Compute residual $\boldsymbol{r} = {\mL_d} \bu - \ba$
        \State Solve linear system $ J^TJ  \boldsymbol{\delta}  = J^T \boldsymbol{r} $ for $\boldsymbol{\delta}$
        \State $\bu_{new} = \bu + \boldsymbol{\delta}$ \\
    \Return{$\bu_{new}$}
	\end{algorithmic} 
\end{algorithm}

\newpage

\bibliographystyle{plain}
\bibliography{Reference}

\begin{thebibliography}{10}

\bibitem{apte2023diffusion}
Rucha Apte, Sheel Nidhan, Rishikesh Ranade, and Jay Pathak.
\newblock Diffusion model based data generation for partial differential equations.
\newblock {\em arXiv preprint arXiv:2306.11075}, 2023.

\bibitem{arisaka2023principled}
Sohei Arisaka and Qianxiao Li.
\newblock Principled acceleration of iterative numerical methods using machine learning.
\newblock {\em Proceedings of the 40th International Conference on Machine Learning}, 2023.

\bibitem{azulay2022multigrid}
Yael Azulay and Eran Treister.
\newblock Multigrid-augmented deep learning preconditioners for the helmholtz equation.
\newblock {\em SIAM Journal on Scientific Computing}, (0):S127--S151, 2022.

\bibitem{bano2020configuration}
Jorge Ba{\~n}o-Medina, Rodrigo Manzanas, and Jos{\'e}~Manuel Guti{\'e}rrez.
\newblock Configuration and intercomparison of deep learning neural models for statistical downscaling.
\newblock {\em Geoscientific Model Development}, 13(4):2109--2124, 2020.

\bibitem{cai2021physics}
Shengze Cai, Zhiping Mao, Zhicheng Wang, Minglang Yin, and George~Em Karniadakis.
\newblock Physics-informed neural networks (pinns) for fluid mechanics: A review.
\newblock {\em Acta Mechanica Sinica}, 37(12):1727--1738, 2021.

\bibitem{chen2022meta}
Yuyan Chen, Bin Dong, and Jinchao Xu.
\newblock Meta-mgnet: Meta multigrid networks for solving parameterized partial differential equations.
\newblock {\em Journal of computational physics}, 455:110996, 2022.

\bibitem{dhariwal2021diffusion}
Prafulla Dhariwal and Alexander Nichol.
\newblock Diffusion models beat gans on image synthesis.
\newblock {\em Advances in neural information processing systems}, 34:8780--8794, 2021.

\bibitem{farimani2017deep}
Amir~Barati Farimani, Joseph Gomes, and Vijay~S Pande.
\newblock Deep learning the physics of transport phenomena.
\newblock {\em arXiv preprint arXiv:1709.02432}, 2017.

\bibitem{goswami2022physics}
Somdatta Goswami, Aniruddha Bora, Yue Yu, and George~Em Karniadakis.
\newblock Physics-informed neural operators.
\newblock {\em arXiv preprint arXiv:2207.05748}, 2022.

\bibitem{goswami2023physics}
Somdatta Goswami, Aniruddha Bora, Yue Yu, and George~Em Karniadakis.
\newblock Physics-informed deep neural operator networks.
\newblock In {\em Machine Learning in Modeling and Simulation: Methods and Applications}, pages 219--254. Springer, 2023.

\bibitem{groenke2020climalign}
Brian Groenke, Luke Madaus, and Claire Monteleoni.
\newblock Climalign: Unsupervised statistical downscaling of climate variables via normalizing flows.
\newblock In {\em Proceedings of the 10th International Conference on Climate Informatics}, pages 60--66, 2020.

\bibitem{han2018solving}
Jiequn Han, Arnulf Jentzen, and Weinan E.
\newblock Solving high-dimensional partial differential equations using deep learning.
\newblock {\em Proceedings of the National Academy of Sciences}, 115(34):8505--8510, 2018.

\bibitem{han2017deep}
Jiequn Han, Arnulf Jentzen, et~al.
\newblock Deep learning-based numerical methods for high-dimensional parabolic partial differential equations and backward stochastic differential equations.
\newblock {\em Communications in mathematics and statistics}, 5(4):349--380, 2017.

\bibitem{ho2020denoising}
Jonathan Ho, Ajay Jain, and Pieter Abbeel.
\newblock Denoising diffusion probabilistic models.
\newblock {\em Advances in Neural Information Processing Systems}, 33:6840--6851, 2020.

\bibitem{ho2022cascaded}
Jonathan Ho, Chitwan Saharia, William Chan, David~J Fleet, Mohammad Norouzi, and Tim Salimans.
\newblock Cascaded diffusion models for high fidelity image generation.
\newblock {\em The Journal of Machine Learning Research}, 23(1):2249--2281, 2022.

\bibitem{hsieh2018learning}
Jun-Ting Hsieh, Shengjia Zhao, Stephan Eismann, Lucia Mirabella, and Stefano Ermon.
\newblock Learning neural pde solvers with convergence guarantees.
\newblock In {\em International Conference on Learning Representations}, 2018.

\bibitem{jiang2023efficient}
Peishi Jiang, Zhao Yang, Jiali Wang, Chenfu Huang, Pengfei Xue, TC~Chakraborty, Xingyuan Chen, and Yun Qian.
\newblock Efficient super-resolution of near-surface climate modeling using the fourier neural operator.
\newblock {\em Journal of Advances in Modeling Earth Systems}, 15(7):e2023MS003800, 2023.

\bibitem{jin2023asymptotic}
Shi Jin, Zheng Ma, and Keke Wu.
\newblock Asymptotic-preserving neural networks for multiscale time-dependent linear transport equations.
\newblock {\em Journal of Scientific Computing}, 94(3):57, 2023.

\bibitem{joshi2019generative}
Ameya Joshi, Viraj Shah, Sambuddha Ghosal, Balaji Pokuri, Soumik Sarkar, Baskar Ganapathysubramanian, and Chinmay Hegde.
\newblock Generative models for solving nonlinear partial differential equations.
\newblock In {\em Proc. of NeurIPS Workshop on ML for Physics}, 2019.

\bibitem{kharazmi2019variational}
Ehsan Kharazmi, Zhongqiang Zhang, and George~Em Karniadakis.
\newblock Variational physics-informed neural networks for solving partial differential equations.
\newblock {\em arXiv preprint arXiv:1912.00873}, 2019.

\bibitem{kharazmi2021hp}
Ehsan Kharazmi, Zhongqiang Zhang, and George~Em Karniadakis.
\newblock hp-vpinns: Variational physics-informed neural networks with domain decomposition.
\newblock {\em Computer Methods in Applied Mechanics and Engineering}, 374:113547, 2021.

\bibitem{kovachki2021neural}
Nikola Kovachki, Zongyi Li, Burigede Liu, Kamyar Azizzadenesheli, Kaushik Bhattacharya, Andrew Stuart, and Anima Anandkumar.
\newblock Neural operator: Learning maps between function spaces.
\newblock {\em arXiv preprint arXiv:2108.08481}, 2021.

\bibitem{leinonen2020stochastic}
Jussi Leinonen, Daniele Nerini, and Alexis Berne.
\newblock Stochastic super-resolution for downscaling time-evolving atmospheric fields with a generative adversarial network.
\newblock {\em IEEE Transactions on Geoscience and Remote Sensing}, 59(9):7211--7223, 2020.

\bibitem{li2020fourier}
Zongyi Li, Nikola Kovachki, Kamyar Azizzadenesheli, Burigede Liu, Kaushik Bhattacharya, Andrew Stuart, and Anima Anandkumar.
\newblock Fourier neural operator for parametric partial differential equations.
\newblock {\em arXiv preprint arXiv:2010.08895}, 2020.

\bibitem{li2021physics}
Zongyi Li, Hongkai Zheng, Nikola Kovachki, David Jin, Haoxuan Chen, Burigede Liu, Kamyar Azizzadenesheli, and Anima Anandkumar.
\newblock Physics-informed neural operator for learning partial differential equations.
\newblock {\em arXiv preprint arXiv:2111.03794}, 2021.

\bibitem{lu2022priori}
Jianfeng Lu and Yulong Lu.
\newblock A priori generalization error analysis of two-layer neural networks for solving high dimensional schr{\"o}dinger eigenvalue problems.
\newblock {\em Communications of the American Mathematical Society}, 2(1):1--21, 2022.

\bibitem{lu2019deeponet}
Lu~Lu, Pengzhan Jin, and George~Em Karniadakis.
\newblock Deeponet: Learning nonlinear operators for identifying differential equations based on the universal approximation theorem of operators.
\newblock {\em arXiv preprint arXiv:1910.03193}, 2019.

\bibitem{lu2021priori}
Yulong Lu, Jianfeng Lu, and Min Wang.
\newblock A priori generalization analysis of the deep ritz method for solving high dimensional elliptic partial differential equations.
\newblock In {\em Conference on learning theory}, pages 3196--3241. PMLR, 2021.

\bibitem{lu2022solving}
Yulong Lu, Li~Wang, and Wuzhe Xu.
\newblock Solving multiscale steady radiative transfer equation using neural networks with uniform stability.
\newblock {\em Research in the Mathematical Sciences}, 9(3):45, 2022.

\bibitem{nikolopoulos2024ai}
Stefanos Nikolopoulos, Ioannis Kalogeris, George Stavroulakis, and Vissarion Papadopoulos.
\newblock Ai-enhanced iterative solvers for accelerating the solution of large-scale parametrized systems.
\newblock {\em International Journal for Numerical Methods in Engineering}, 125(2):e7372, 2024.

\bibitem{price2022increasing}
Ilan Price and Stephan Rasp.
\newblock Increasing the accuracy and resolution of precipitation forecasts using deep generative models.
\newblock In {\em International conference on artificial intelligence and statistics}, pages 10555--10571. PMLR, 2022.

\bibitem{raissi2019physics}
Maziar Raissi, Paris Perdikaris, and George~E Karniadakis.
\newblock Physics-informed neural networks: A deep learning framework for solving forward and inverse problems involving nonlinear partial differential equations.
\newblock {\em Journal of Computational physics}, 378:686--707, 2019.

\bibitem{ronneberger2015u}
Olaf Ronneberger, Philipp Fischer, and Thomas Brox.
\newblock U-net: Convolutional networks for biomedical image segmentation.
\newblock In {\em Medical Image Computing and Computer-Assisted Intervention--MICCAI 2015: 18th International Conference, Munich, Germany, October 5-9, 2015, Proceedings, Part III 18}, pages 234--241. Springer, 2015.

\bibitem{sachindra2018statistical}
DA~Sachindra, Khandakar Ahmed, Md~Mamunur Rashid, S~Shahid, and BJC Perera.
\newblock Statistical downscaling of precipitation using machine learning techniques.
\newblock {\em Atmospheric research}, 212:240--258, 2018.

\bibitem{shu2023physics}
Dule Shu, Zijie Li, and Amir~Barati Farimani.
\newblock A physics-informed diffusion model for high-fidelity flow field reconstruction.
\newblock {\em Journal of Computational Physics}, 478:111972, 2023.

\bibitem{song2020denoising}
Jiaming Song, Chenlin Meng, and Stefano Ermon.
\newblock Denoising diffusion implicit models.
\newblock {\em arXiv preprint arXiv:2010.02502}, 2020.

\bibitem{song2020score}
Yang Song, Jascha Sohl-Dickstein, Diederik~P Kingma, Abhishek Kumar, Stefano Ermon, and Ben Poole.
\newblock Score-based generative modeling through stochastic differential equations.
\newblock {\em arXiv preprint arXiv:2011.13456}, 2020.

\bibitem{um2020solver}
Kiwon Um, Robert Brand, Yun~Raymond Fei, Philipp Holl, and Nils Thuerey.
\newblock Solver-in-the-loop: Learning from differentiable physics to interact with iterative pde-solvers.
\newblock {\em Advances in Neural Information Processing Systems}, 33:6111--6122, 2020.

\bibitem{vandal2019intercomparison}
Thomas Vandal, Evan Kodra, and Auroop~R Ganguly.
\newblock Intercomparison of machine learning methods for statistical downscaling: the case of daily and extreme precipitation.
\newblock {\em Theoretical and Applied Climatology}, 137:557--570, 2019.

\bibitem{vandal2017deepsd}
Thomas Vandal, Evan Kodra, Sangram Ganguly, Andrew Michaelis, Ramakrishna Nemani, and Auroop~R Ganguly.
\newblock Deepsd: Generating high resolution climate change projections through single image super-resolution.
\newblock In {\em Proceedings of the 23rd acm sigkdd international conference on knowledge discovery and data mining}, pages 1663--1672, 2017.

\bibitem{wang2023long}
Sifan Wang and Paris Perdikaris.
\newblock Long-time integration of parametric evolution equations with physics-informed deeponets.
\newblock {\em Journal of Computational Physics}, 475:111855, 2023.

\bibitem{wang2023expert}
Sifan Wang, Shyam Sankaran, Hanwen Wang, and Paris Perdikaris.
\newblock An expert's guide to training physics-informed neural networks.
\newblock {\em arXiv preprint arXiv:2308.08468}, 2023.

\bibitem{wang2021learning}
Sifan Wang, Hanwen Wang, and Paris Perdikaris.
\newblock Learning the solution operator of parametric partial differential equations with physics-informed deeponets.
\newblock {\em Science advances}, 7(40):eabi8605, 2021.

\bibitem{wei2023super}
Min Wei and Xuesong Zhang.
\newblock Super-resolution neural operator.
\newblock In {\em Proceedings of the IEEE/CVF Conference on Computer Vision and Pattern Recognition}, pages 18247--18256, 2023.

\bibitem{weinan2018deep}
E~Weinan and Bing Yu.
\newblock The deep ritz method: A deep learning-based numerical algorithm for solving variational problems.
\newblock {\em Communications in Mathematics and Statistics}, 6(1):1--12, 2018.

\bibitem{wilby1998statistical}
Robert~L Wilby, TML Wigley, D~Conway, PD~Jones, BC~Hewitson, J~Main, and DS~Wilks.
\newblock Statistical downscaling of general circulation model output: A comparison of methods.
\newblock {\em Water resources research}, 34(11):2995--3008, 1998.

\bibitem{yang2023denoising}
Gefan Yang and Stefan Sommer.
\newblock A denoising diffusion model for fluid field prediction.
\newblock {\em arXiv e-prints}, pages arXiv--2301, 2023.

\bibitem{yang2023fourier}
Qidong Yang, Alex Hernandez-Garcia, Paula Harder, Venkatesh Ramesh, Prasanna Sattegeri, Daniela Szwarcman, Campbell~D Watson, and David Rolnick.
\newblock Fourier neural operators for arbitrary resolution climate data downscaling.
\newblock {\em arXiv preprint arXiv:2305.14452}, 2023.

\bibitem{yu2022gradient}
Jeremy Yu, Lu~Lu, Xuhui Meng, and George~Em Karniadakis.
\newblock Gradient-enhanced physics-informed neural networks for forward and inverse pde problems.
\newblock {\em Computer Methods in Applied Mechanics and Engineering}, 393:114823, 2022.

\bibitem{zang2020weak}
Yaohua Zang, Gang Bao, Xiaojing Ye, and Haomin Zhou.
\newblock Weak adversarial networks for high-dimensional partial differential equations.
\newblock {\em Journal of Computational Physics}, 411:109409, 2020.

\end{thebibliography}

\end{document}